\documentclass[12pt]{article}
\usepackage{amsthm,amsmath,amssymb,amscd,amsfonts, bm}
\begin{document}
\centerline{\bf Splitting of Unstable 2-Bundles}
\centerline{\bf Over the Complex Projective 6-Space}

\bigbreak\centerline{\it Dedicated to Thomas Peternell}

\bigbreak
\centerline{Yum-Tong Siu\ %
\footnote{Partially supported by grant DMS-1001416 from the National Science Foundation.}
}

\bigbreak

\bigbreak\noindent{\bf\S1. Introduction.}  In this article we prove the following theorem.

\bigbreak\noindent(1.1) {\it Main Theorem.} Let $V$ be an unstable
holomorphic 2-bundle over ${\mathbb P}_n$ for some $n\geq 6$. Then
$V$ is a direct sum of two holomorphic line bundles.

\bigbreak The statement of our Main Theorem with the weaker dimension assumption of $n\geq 4$ was already given in the paper of Grauert-Schneider [Grauert-Schneider1977] in 1977, but some points of the proof seem not to have yet been clearly worked out there (see the review of Wolf Barth [Barth1979]). Schneider in his survey paper [Schneider1987, p.104, lines 1-4] ten years later posed the case of $n\geq 5$ as a conjecture.  Since 1977 there has not yet been any really substantial progress toward the resolution of the conjecture.

\medbreak Background material of the theory of vector bundles over complex projective spaces and the significance of the conjecture of decomposability of unstable $2$-bundles on ${\mathbb P}_n$ for $n\geq 4$ can be found in [Okonek-Schneider-Spindler2011, vandeVen1980, Schneider1987].

\bigbreak\noindent(1.2) {\it New Ingredients in the Method of Proof.}  Compared to the approach of the 1977 paper of Grauert-Schneider [Grauert-Schneider1977], our method of proof uses the following three new ingredients.

\medbreak\noindent(i) The 1977 paper of Grauert-Schneider uses the positivity of the tangent bundle of ${\mathbb P}_n$ and vanishing theorems for the unreduced codimension 2 subspace defined from the instability of the $2$-bundle $V$.  Complications arise from the unreduced structure of the codimension 2 subspace.  Here, instead of using the positivity of the tangent bundle of ${\mathbb P}_n$, we use holomorphic vector fields on ${\mathbb P}_n$ to directly construct holomorphic sections of the given bundle on the reduction of the codimension 2 subspace.  Our arguments avoid the complications in the application of vanishing theorems which result from the unreduced structure of the codimension 2 subspace. This ingredient is used from \S2 on.

\medbreak\noindent(ii) We use Mathias Peternell's singular subvariety version [Peternell1983] of the theorem of Barth-Lefschetz Theorem [Barth1970, Barth-Larsen1972, Larsen1973, Barth1975, Schneider1975, Sommese1979, Sommese1982, Schneider-Zintl1993] to extend a line bundle on a branch of the reduction of the codimension 2 subspace to ${\mathbb P}_n$.  This step requires $n\geq 6$.  An adaptation of our method to $n=5$ would require the analysis of the critical dimension for our special situation and not just applying the ready-made result of the theorem of Barth-Lefschetz-Peternell.  This ingredient is used from \S3 on.

\medbreak\noindent(iii) We introduce a (possibly singular) branched cover of the blow-up of ${\mathbb P}_n$ along a complex line to handle the problem of the normal sheaf of a branch of the reduction of the codimension 2 subspace not being locally free.  Instead of working with a branch of the reduction of the codimension 2 subspace whose normal sheaf is not locally free, the (possibly singular) branched cover of the blow-up of ${\mathbb P}_n$ along a complex line enables us to work with the locally free normal sheaf of a singular codimension 2 subvariety in it.  The blow-up of ${\mathbb P}_n$ along a complex line occurs in the context of the light-source projection from the complement of the complex line in ${\mathbb P}_n$ to a linear subspace ${\mathbb P}_{n-2}$ in ${\mathbb P}_n$ which is disjoint from the complex line.  This ingredient is used in \S4.

\bigbreak\noindent(1.3) {\it Organization of the Arguments of the Proof.}  To better present and highlight the roles of the different ingredients in the proof, we give the arguments of the proof first in the following three special cases. (i) The codimension 2 subspace $Z$ defined from the instability of the $2$-bundle $V$ is reduced at some point.  (ii) The two generic vanishing orders of the codimension 2 subspace $Z$ are equal at some point. (iii) One branch of the reduction of the codimension 2 subspace is regular.  For the special cases (i) and (ii) it suffices to assume the complex dimension $n$ of the complex projective space ${\mathbb P}_n$ to be at least $3$.  The proof of the special case (iii) already requires $n$ to be $\geq 6$.  After the presentation of the three special cases, the proof of the general case is given with the emphasis on how to handle the singularities of a branch of the reduction of the codimension 2 subspace.  The discussion of various special cases enables us to see in isolation the role of each new ingredient of the method of the proof.

\medbreak\noindent(1.3.1) {\it Notations.}  The notations ${\mathbb Z}$ and ${\mathbb N}$ mean respectively the set of all integers and the set of all positive integers.  For a complex space $X$ the notation ${\mathcal O}_X$ means the structure sheaf of $X$.  An unreduced complex space $(X,{\mathcal O}_X)$ means that nonzero nilpotent elements are allowed in the structure sheaf ${\mathcal O}_X$.  For an unreduced complex space $Z$ the notation $Z_{\rm red}$ means the reduction of $Z$ (which is obtained by replacing the structure sheaf by its quotient by the subsheaf of all nilpotent elements).  For a coherent sheaf ${\mathcal F}$ on a complex space $X$, the notation $\Gamma\left(X, {\mathcal F}\right)$ means the vector space of all global sections of ${\mathcal F}$ over $X$.  For a holomorphic vector bundle $E$ on a complex space $X$ the notation $E^*$ means the dual bundle of $E$, the notation ${\mathcal O}_X(E)$ means the sheaf of germs of holomorphic sections of $E$, and the notation $\Gamma(X,E)$ means $\Gamma(X,{\mathcal O}_X(E))$.  A reduced complex space $W$ is called a {\it branched-cover} of a reduced irreducible complex space $W$ means that there is a surjective holomorphic map $\pi$ from $\hat W$ to $W$ with finite fibers which is a covering map outside a proper subvariety of $W$.  The map $\pi$ is called a {\it branched-cover map}.  For $\ell\in{\mathbb Z}$ the notation ${\mathcal O}_{{\mathbb P}_n}(\ell)$ means the $\ell$-th tensor power of the hyperplane section line bundle on ${\mathbb P}_n$.  For a coherent sheaf ${\mathcal F}$ on ${\mathbb P}_n$ and $\ell\in{\mathbb Z}$ the notation ${\mathcal F}(\ell)$ means the tensor product of
${\mathcal F}$ and ${\mathcal O}_{{\mathbb P}_n}(\ell)$.  For a complex submanifold $A$ in a complex manifold $B$, the notation $N_{A,B}$ means the normal bundle of $A$ in $B$.  

\bigbreak\noindent(1.4) {\it Codimension 2 Subspace from Instability of $2$-Bundle.}
Assume that the given holomorphic unstable $2$-bundle $V$ over the complex projective space
${\mathbb P}_n$ is not a direct sum of two holomorphic
line bundles over ${\mathbb P}_n$.  We are going to derive a contradiction from the assumption of $n\geq 6$ and in some special cases from the weaker assumption of $n\geq 3$.  We first introduce the codimension 2 subspace $Z$ in ${\mathbb P}_n$ defined from the instability of the $2$-bundle $V$.  In general, the subspace $Z$ is unreduced.  Since $V$ is unstable, there
exists a sheaf-monomorphism $\Phi$ from ${\mathcal O}_{{\mathbb
P}_n}(\ell)$ to ${\mathcal O}_{{\mathbb P}_n}(V)$ such that

\medbreak\noindent(a) the cokernel of $\Phi$ is locally free outside a
subvariety $Z$ of complex codimension $2$ in ${\mathbb
P}_n$, and

\medbreak\noindent(b) the first Chern class of the determinant line bundle of
$V$ is no greater than $2\ell$.

\medbreak\noindent By replacing $V$ by $V(-\ell)$, we can assume without loss of
generality that $\ell=0$ so that we have a sheaf-monomorphism
$\Phi:{\mathcal O}_{{\mathbb P}_n}\to{\mathcal O}_{{\mathbb
P}_n}(V)$ (which satisfies the condition (a)) and the first Chern class of $V$ is no greater than $0$.  Let
$b$ be the nonnegative integer such that the first Chern class of $V$ is $-b$.

\medbreak Let $s\in\Gamma\left({\mathbb P}_n,V\right)$ such that
the image of $\Phi$ is generated by $s$ over ${\mathcal
O}_{{\mathbb P}_n}$.  We can cover ${\mathbb P}_n$ by a finite
number of open subsets $U_j$ ($1\leq j\leq J$) such that

\medbreak\noindent (i) the transition function for the vector
bundle $V$ from $U_k$ to $U_j$ is the $2\times 2$ matrix
$$
\left(\begin{matrix}a_{jk}&b_{jk}\cr
c_{jk}&d_{jk}\end{matrix}\right)
$$
of holomorphic functions on $U_j\cap U_k$,

\medbreak\noindent (ii) $s|_{U_j}$ is represented by the column
$2$-vector
$$
\left(\begin{matrix}\text{$\mathfrak f_j$}\cr {\mathfrak
g_j}\end{matrix}\right)
$$
of holomorphic functions on $U_j$, and
$$
\left(\begin{matrix}{\mathfrak f_j}\cr {\mathfrak
g_j}\end{matrix}\right) =\left(\begin{matrix}a_{jk}&b_{jk}\cr
c_{jk}&d_{jk}\end{matrix}\right)\left(\begin{matrix}{\mathfrak
f_k}\cr {\mathfrak g_k}\end{matrix}\right)\quad{\rm on}\ \ U_j\cap U_k. \leqno{(1.4.1)}
$$

\bigbreak\noindent We define the complex subspace $Z$ of ${\mathbb
P}_n$ endowed with the complex structure ${\mathcal O}_Z$ so that
$${\mathcal O}_Z={\mathcal O}_{{\mathbb P}_n}\left/\left({\mathcal
O}_{{\mathbb P}_n}{\mathfrak f_j}+{\mathcal O}_{{\mathbb
P}_n}{\mathfrak g_j}\right)\right.
$$
on $U_j\cap Z$.

\bigbreak

\bigbreak\noindent{\bf\S2. Technique of Vector Fields.}

\bigbreak\noindent(2.1) {\it Special Case of the Codimension 2 Subspace Reduced at Some Point}. If the structure sheaf of $\left(Z,{\mathcal O}_Z\right)$ is
reduced at some point $P_0$ of $Z$ (which we can assume without loss of generality to be
a regular point of $Z$), then we can easily derive a
contradiction in the following way even for the case of $n\geq 3$.  We take any generic holomorphic vector field
$\xi$ on ${\mathbb P}_n$ (which vanishes on some hypersurface $H$ of
${\mathbb P}_n$) and apply $\xi$ to both sides of the equation $(1.4.1)$ to get
$$
\left(\begin{matrix}\xi\left({\mathfrak f_j}\right)\cr
\xi\left({\mathfrak g_j}\right)\end{matrix}\right)
=\left(\begin{matrix}a_{jk}&b_{jk}\cr
c_{jk}&d_{jk}\end{matrix}\right)\left(\begin{matrix}\xi\left({\mathfrak
f_k}\right)\cr \xi\left({\mathfrak
g_k}\right)\end{matrix}\right)\leqno{(2.1.1)}
$$
modulo ${\mathcal O}_{{\mathbb P}_n}{\mathfrak f_j}+{\mathcal
O}_{{\mathbb P}_n}{\mathfrak g_j}$ on $U_j\cap U_k$.  We take another generic holomorphic vector field $\eta$ on ${\mathbb P}_n$
(which vanishes on some hypersurface $H$ of ${\mathbb P}_n$) to get likewise
$$
\left(\begin{matrix}\eta\left({\mathfrak f_j}\right)\cr
\eta\left({\mathfrak g_j}\right)\end{matrix}\right)
=\left(\begin{matrix}a_{jk}&b_{jk}\cr
c_{jk}&d_{jk}\end{matrix}\right)\left(\begin{matrix}\eta\left({\mathfrak
f_k}\right)\cr \eta\left({\mathfrak
g_k}\right)\end{matrix}\right)\leqno{(2.1.2)}
$$
modulo ${\mathcal O}_{{\mathbb P}_n}{\mathfrak f_j}+{\mathcal
O}_{{\mathbb P}_n}{\mathfrak g_j}$ on $U_j\cap U_k$.  Putting
together $(2.1.1)$ and $(2.1.2)$, we get
$$
\left(\begin{matrix}\xi\left({\mathfrak
f_j}\right)&\eta\left({\mathfrak f_j}\right)\cr
\xi\left({\mathfrak g_j}\right)&\eta\left({\mathfrak
g_j}\right)\end{matrix}\right)
=\left(\begin{matrix}a_{jk}&b_{jk}\cr
c_{jk}&d_{jk}\end{matrix}\right)\left(\begin{matrix}\xi\left({\mathfrak
f_k}\right)&\eta\left({\mathfrak f_k}\right)\cr
\xi\left({\mathfrak g_j}\right)&\eta\left({\mathfrak
g_k}\right)\end{matrix}\right)\leqno{(2.1.3)}
$$
modulo ${\mathcal O}_{{\mathbb P}_n}{\mathfrak f_j}+{\mathcal
O}_{{\mathbb P}_n}{\mathfrak g_j}$ on $U_j\cap U_k$.  By taking the determinants of both sides of (2.1.3), we obtain
$$
\det\left(\begin{matrix}\xi\left({\mathfrak
f_j}\right)&\eta\left({\mathfrak f_j}\right)\cr
\xi\left({\mathfrak g_j}\right)&\eta\left({\mathfrak
g_j}\right)\end{matrix}\right)
=\det\left(\begin{matrix}a_{jk}&b_{jk}\cr
c_{jk}&d_{jk}\end{matrix}\right)\
\det\left(\begin{matrix}\xi\left({\mathfrak
f_k}\right)&\eta\left({\mathfrak f_k}\right)\cr
\xi\left({\mathfrak g_j}\right)&\eta\left({\mathfrak
g_k}\right)\end{matrix}\right)
$$
modulo ${\mathcal O}_{{\mathbb P}_n}{\mathfrak f_j}+{\mathcal
O}_{{\mathbb P}_n}{\mathfrak g_j}$ on $U_j\cap U_k$.  With generic choices of the holomorphic vector fields $\xi$ and $\eta$ on ${\mathbb P}_n$, we have a
non identically zero element of $\Gamma\left(Z,{\mathcal
O}_Z\left(\det V\right)\right)$ which vanishes on $H\cap Z$. In particular, on the branch $Z^0_{\rm red}$ of the reduction
$Z_{\rm red}$ of $Z$ which contains the point $P_0$ we have a non trivial global
holomorphic section $\sigma$ of $\det V$.  This contradicts the fact
that the first Chern class of $V$ is nonpositive, because the non identically zero holomorphic section $\sigma$ of $\det V$ on $Z^0_{\rm red}$ vanishes on the ample divisor $H\cap Z^0_{\rm red}$ of $Z^0_{\rm red}$.  The dimension assumption $n\geq 3$ is used
to guarantee that the complex dimension of the irreducible subvariety $Z^0_{\rm red}$ is at least $1$.

\bigbreak\noindent(2.2) {\it Two Vanishing Orders of Defining Ideal Sheaf of Codimension $2$ Subspace.}  At a generic point $P_0$ 
of $Z$ (in some $U_k$) we can take a local surface $\Omega$ which is transversal
to (the reduction of) $Z$ at $P_0$ and is biholomorphic to the bidisk $\Delta^2\subset{\mathbb
C}^2$ whose coordinates are $\left(\zeta_1,\zeta_2\right)$ with $P_0$ corresponding to the origin $(\zeta_1,\zeta_2)=(0,0)$.  Consider the two holomorphic functions
$\phi\left(\zeta_1,\zeta_2\right)$ and
$\psi\left(\zeta_1,\zeta_2\right)$ which are respectively the
restrictions of ${\mathfrak f}_k$ and ${\mathfrak g}_k$ to
$\Omega$.  Consider their expansions into homogeneous parts
$$
\displaylines{\phi\left(\zeta_1,\zeta_2\right)=\sum_{\nu=\kappa_1}^\infty\phi_\nu\left(\zeta_1,\zeta_2\right),
\cr
\psi\left(\zeta_1,\zeta_2\right)=\sum_{\nu=\kappa_2}^\infty\psi_\nu\left(\zeta_1,\zeta_2\right).
\cr}
$$
where $\phi_{\kappa_1}$ and $\psi_{\kappa_2}$ are both non
identically zero and not equal to nonzero constants times each
other.  If $\psi_{\kappa_2}=c\,\phi_{\kappa_1}$ for $c\in{\mathbb C}$
at a generic point of $Z\cap U_k$, we change
local basis of $V$ on $U_k$ by
$$
\left\{\begin{matrix}{\mathfrak f}_k\to {\mathfrak f}_k\hfill\cr
{\mathfrak g}_k\to {\mathfrak g}_k- c\,{\mathfrak
f}_k\cr\end{matrix}\right.$$
to make sure that $\phi_{\kappa_1}$ and
$\psi_{\kappa_2}$ are not equal to nonzero constants times each
other when $\kappa_1=\kappa_2$.
We call $\kappa_1$ and $\kappa_2$
the {\it two generic vanishing orders} for
the defining sheaf for $Z$ (at some generic point of $Z$ in $U_k$).  We agree to use the labeling so that $\kappa_1\leq\kappa_2$.

\bigbreak\noindent(2.3) {\it Special Case of Two Equal Generic Vanishing
Orders.} The above argument of using global vector fields of ${\mathbb P}_n$ to conclude
the nonexistence of a nonsplitting unstable $2$-bundle on ${\mathbb P}_n$ with $n\geq 3$ in the special case of (2.1) works also in the special case where the two generic vanishing orders
of the defining ideal sheaf of $Z$ are equal at some point of $Z$.  Of course, the case of  (2.1) where the codimension 2 subspace $Z$ is reduced at some generic point simply means that both generic vanishing orders are equal to $1$ at that generic point.

\medbreak  Suppose the generic vanishing
orders of the two local defining functions for the ideal sheaf of $Z$ are both equal to $m\geq 2$
at some generic point $P_0$ of $Z$. Then we can choose
global holomorphic vector fields $\xi_1,\cdots,\xi_{m-1}$ and
$\eta_1,\cdots,\eta_{m-1}$ on ${\mathbb P}_n$ vanishing on some hyperplane $H$ of ${\mathbb P}_n$ so that
$$
\begin{aligned}&\det\left(\begin{matrix}\xi_1\cdots\xi_{m-1}\left({\mathfrak
f_j}\right)&\eta_1\cdots\eta_{m-1}\left({\mathfrak f_j}\right)\cr
\xi_1\cdots\xi_{m-1}\left({\mathfrak
g_j}\right)&\eta_1\cdots\eta_{m-1}\left({\mathfrak
g_j}\right)\end{matrix}\right)\cr
&=
\det\left(\begin{matrix}a_{jk}&b_{jk}\cr
c_{jk}&d_{jk}\end{matrix}\right)
\det\left(\begin{matrix}\xi_1\cdots\xi_{m-1}\left({\mathfrak
f_k}\right)&\eta_1\cdots\eta_{m-1}\left({\mathfrak f_k}\right)\cr
\xi_1\cdots\xi_{m-1}\left({\mathfrak
g_k}\right)&\eta_1\cdots\eta_{m-1}\left({\mathfrak
g_k}\right)\end{matrix}\right)\cr
\end{aligned}
$$
on the reduction $Z_{\rm red}$ of $Z$.  The column $2$-vector
$$
\det\left(\begin{matrix}\xi_1\cdots\xi_{m-1}\left({\mathfrak
f_j}\right)&\eta_1\cdots\eta_{m-1}\left({\mathfrak f_j}\right)\cr
\xi_1\cdots\xi_{m-1}\left({\mathfrak
g_j}\right)&\eta_1\cdots\eta_{m-1}\left({\mathfrak
g_j}\right)\end{matrix}\right)
$$
of holomorphic functions on $U_j$ for $1\leq j\leq J$ gives us a nontrivial global
holomorphic section $\sigma$ of $\det V$ on some branch $Z_{\rm red}^0$
of $Z_{\rm red}$.  This 
contradicts the fact that the first Chern class of $V$ is
nonpositive, because $\sigma$ vanishes on the ample divisor $H\cap Z^0_{\rm red}$ of $Z^0_{\rm red}$. Again for this present situation of two equal
generic vanishing orders the dimension assumption $n\geq 3$ is used
to guarantee that the complex dimension of the irreducible subvariety $Z^0_{\rm red}$ is at least $1$.

\bigbreak\noindent(2.4) {\it Construction of Section for Case of Two Unequal Generic Vanishing Orders.}  When the two generic vanishing orders of the defining ideal sheaf for the subspace $Z$ are unequal at a generic point of any branch of the reduction $Z_{\rm red}$ of $Z$, we can still let $m$ be the smallest of the two generic vanishing orders at the generic point of any branch of ${\rm red}\,Z$ and choose
global holomorphic vector fields $\xi_1,\cdots,\xi_{m-1}$ on ${\mathbb P}_n$ vanishing on some hyperplane $H$ of ${\mathbb P}_n$ so that
$$
\left(\begin{matrix}\xi_1\cdots\xi_{m-1}\left({\mathfrak
f_j}\right)\cr
\xi_1\cdots\xi_{m-1}\left({\mathfrak
g_j}\right)\cr\end{matrix}\right)=
\left(\begin{matrix}a_{jk}&b_{jk}\cr
c_{jk}&d_{jk}\end{matrix}\right)
\left(\begin{matrix}\xi_1\cdots\xi_{m-1}\left({\mathfrak
f_k}\right)\cr
\xi_1\cdots\xi_{m-1}\left({\mathfrak
g_k}\right)\cr\end{matrix}\right)
$$
on $U_j\cap U_k\cap Z_{\rm red}$ yields a non identically zero element $\sigma\in\Gamma\left(Z_{\rm red}, V\right)$.  Of course, for this general case we are unable to use the same method to produce another holomorphic section $\sigma^\prime$ of $V$ over $Z_{\rm red}$ such that $\sigma\wedge\sigma^\prime$ is a non identically zero holomorphic section of $\det V$ over $Z_{\rm red}$.  It is for this reason why we need to introduce the technique involving the Barth-Lefschetz-Peternell theorem which requires the stronger dimension assumption $n\geq 6$.

\bigbreak

\bigbreak\noindent{\bf\S3. Line Bundle Extension by Theorem of Barth-Lefschetz-Peternell.}

\bigbreak\noindent(3.1) {\it Cartier Divisor for Case of Two Unequal Generic Vanishing Orders.}  Assume that the two generic vanishing orders are unequal with $m$ being the smaller one.  Let $Z_{\rm red}^0$ be a branch of the reduction $Z_{\rm red}$ of $Z$.  Let $\sigma^0\in\Gamma\left(Z_{\rm red}^0, V\right)$ be defined by $\sigma\in\Gamma\left(Z_{\rm red}, V\right)$ constructed in (2.4), where the holomorphic vector fields $\xi_1,\cdots,\xi_{m-1}$ on ${\mathbb P}_n$ used in the construction of $\sigma$ in (2.4) are assumed to be so chosen that $\sigma^0$ is not identically zero.  Since $m$ is the smaller of the two unequal generic vanishing orders, with the convention of $\kappa_1\leq\kappa_2$ of (2.2), only the first entry $\xi_1\cdots\xi_{m-1}\left({\mathfrak
f_j}\right)$ in the column $2$-vector
$$
\left(\begin{matrix}\xi_1\cdots\xi_{m-1}\left({\mathfrak
f_j}\right)\cr
\xi_1\cdots\xi_{m-1}\left({\mathfrak
g_j}\right)\cr\end{matrix}\right)
$$
representing $\sigma^0$ is not identically zero on $Z_{\rm red}^0\cap U_j$.  The other entry $\xi_1\cdots\xi_{m-1}\left({\mathfrak
g_j}\right)$ is identically zero on $Z_{\rm red}^0\cap U_j$.  Thus the zero-set of $\sigma^0$ is of pure codimension $1$ in $Z_{\rm red}^0$ and is given on $U_j$ by the Cartier divisor defined by  $\xi_1\cdots\xi_{m-1}\left({\mathfrak
f_j}\right)$, which we also know must contain $H\cap Z_{\rm red}^0$, where $H$ is the hypersurface of ${\mathbb P}_n$ where the holomorphic vector field $\xi_1$ vanishes.
Let $W$ be the Cartier divisor in $Z_{\rm red}^0$ defined by $\sigma^0\in\Gamma\left(Z_{\rm red}^0, V\right)$.  Denote by $[W]$ the holomorphic line bundle on $Z_{\rm red}^0$ defined by the Cartier divisor $W$.

\bigbreak\noindent(3.2) {\it Extension of Line Bundle Defined by Divisor for Special Case of Codimension 2 Subspace Having Regular Branch.}  We now discuss how to handle the general case of two unequal generic vanishing orders
of the defining ideal sheaf for the codimension $2$ subspace $Z$ under the stronger dimension assumption of $n\geq 6$.  In order to better explain the use of Barth-Lefschetz-Peternell theorem and the reason for the later construction of a branched cover over the blow-up of ${\mathbb P}_n$ along a complex line, first we present here the proof for the special case where there is a regular branch $Z_{\rm red}^0$ in the reduction $Z_{\rm red}$ of the unreduced codimension 2 subspace $Z$ defined from the instability of the given indecomposable $2$-bundle $V$.  

\medbreak To simplify notations we denote $Z_{\rm red}^0$ by $A$.    We are going to use the Barth-Lefschetz theorem (before its generalization by Mathias Peternell to singular subvarieties) to show that the line bundle $[W]$ on $A=Z_{\rm red}^0$ from (3.1) is the restriction to $Z_{\rm red}^0$ of the line bundle ${\mathcal O}_{{\mathbb P}_n}(a)$ on ${\mathbb P}_n$ for some $a\in{\mathbb N}$.

\medbreak Since by $2\leq 2\dim_{\mathbb C}A-n$ from $n\geq 6$ and $\dim_{\mathbb C}A=n-2$ the Barth-Lefschetz theorem (in the setting of cohomology groups with coefficients in ${\mathbb Z}$) implies that $H^2\left({\mathbb P}_n, {\mathbb Z}\right)\to H^2(A, {\mathbb Z})$ induced by the inclusion map $A\hookrightarrow{\mathbb P}_n$ is an isomorphism, it follows that there exists some $a\in{\mathbb Z}$ such that the Chern class of the restriction of ${\mathcal O}_{{\mathbb P}_n}(a)$ to $A$ is the Chern class of the line bundle $[W]$ on $A$.  On the compact K\"ahler manifold $A$ the line bundle $[W]\otimes{\mathcal O}_{{\mathbb P}_n}(-a)$ of zero Chern class must be a flat line bundle (in the sense that local fiber coordinates of it can be chosen so that its transition functions are constant functions).   Since $1\leq 2\dim_{\mathbb C}A-n$ from $n\geq 5$ and $\dim_{\mathbb C}A=n-2$, the Barth-Lefschetz theorem (in the setting of homotopy groups) implies that $\pi_1(A)\to\pi({\mathbb P}_n)$ induced by the inclusion map $A\hookrightarrow{\mathbb P}_n$ is an isomorphism and hence $A$ is simply connected.  Since any flat line bundle on a simply connected compact manifold is trivial, it follows that the flat line bundle $[W]\otimes{\mathcal O}_{{\mathbb P}_n}(-a)$ on $A$ is trivial and $[W]$ is the restriction of ${\mathcal O}_{{\mathbb P}_n}(a)$ to $A$.  From $A\cap H\subset W$ in (3.1) we conclude that $a>0$.

\bigbreak\noindent(3.3) {\it Extension of Section from Vanishing of Cohomology.}  Let $t_W$ be the canonical section of $[W]$ so that $t_W$ is a holomorphic section of ${\mathcal O}_{{\mathbb P}_n}(a)$ over $A$, which implies that $\frac{\sigma^0}{t_W}$ as a holomorphic section of $V(-a)$ over $A$ must be nowhere zero on $A$.  Let $\tau=\frac{\sigma^0}{t_W}$.  We have a short exact sequence
$$0\to{\mathcal O}_A\stackrel{\Phi}{\longrightarrow}{\mathcal O}_A(V(-a))\to{\mathcal O}_A(-2a-b)\to 0\leqno{(3.3.1)}$$ on $A$, where the inclusion map $\Phi$ is defined by multiplication by $\tau$, because the Chern class of the determinant line bundle $\det V$ of $V$ is $-b$ so that the Chern class of the determinant line bundle $\det(V(-a))$ of $V$ is $-2a-b$.

\medbreak We are going to extend $\tau$ from $A$ to a global holomorphic section $\tilde\tau$ of $V(-a)$ over all of ${\mathbb P}_n$ in the following way.  Let ${\mathcal I}d_A$ be the full ideal sheaf for the complex submanifold $A$ of ${\mathbb P}_n$.
One key ingredient for the extendibility of $\tau$ is the vanishing
$$H^1\left(A, \left(\left.{\mathcal I}d_A^\ell\right/{\mathcal I}d_A^{\ell+1}\right)(-q)\right)=0\quad{\rm for}\ \ \ell\geq 1\ \ {\rm and}\ \  q\geq 0.
\leqno{(3.3.2)}$$
Since the sheaf $\left.{\mathcal I}d_A^\ell\right/{\mathcal I}d_A^{\ell+1}$ is the sheaf of germs of holomorphic sections of the $\ell$-th symmetric power ${\rm Sym}^\ell\left(\left(N_{A,{\mathbb P}_n}\right)^*\right)$ of the dual $\left(N_{A,{\mathbb P}_n}\right)^*$ of the normal bundle $N_{A,{\mathbb P}_n}$ of the submanifold $A$ of ${\mathbb P}_n$, the vanishing $(3.3.2)$ follows from the result of Schneider-Zintl [Schneider-Zintl1993, p.261, Theorem 4] that for a complex submanifold $Y$ of complex dimension $n_Y$ in ${\mathbb P}_n$,
$$H^p\left(Y, {\rm Sym}^\ell\left(\left(N_{Y,{\mathbb P}_n}\right)^*\right)(-q)\right) = 0
$$
for $q\geq 0$, $\ell\geq 1$ and $p\leq 2n_Y-n$.  In our case $p=1$ and $n_Y=n-2$ so that $1\leq 2(n-2)-n=n-4$ is satisfied with $n\geq 5$.  As noted in the 1993 paper of Schneider-Zintl [Schneider1993, lines 4 and 5 from the bottom of p.259], the result of Schneider-Zintl which we use here is a reproof of Faltings's theorem [Faltings1981, p.145, Korollar 1] in a slightly different formulation by   using the vanishing theorem of Le Potier [LePotier1975].

\medbreak Note that though the normal bundle $N_{Y,{\mathbb P}_n}$ of a complex submanifold $Y$ of ${\mathbb P}_n$ carries positivity from the global holomorphic vector fields of ${\mathbb P}_n$, the vanishing theorem of Kodaira for Nakano positive vector bundles [Nakano1955] cannot be applied to this kind of positivity of $N_{Y,{\mathbb P}_n}$ to give the vanishing of $H^p\left(Y, {\rm Sym}^\ell\left(\left(N_{Y,{\mathbb P}_n}\right)^*\right)(-q)\right)$
for $q\geq 0$, $\ell\geq 1$ and $p\leq 2n_Y-n$ in the result of Schneider-Zintl.

\bigbreak Tensoring the short exact sequence $(3.3.1)$ with ${\rm Sym}^\ell\left(\left(N_{A,{\mathbb P}_n}\right)^*\right)$, we get the short exact sequence

{\footnotesize$$0\to\left.{\mathcal I}d_A^\ell\right/{\mathcal I}d_A^{\ell+1}\to\left(\left.{\mathcal I}d_A^\ell\right/{\mathcal I}d_A^{\ell+1}\right)(V(-a))\to\left(\left.{\mathcal I}d_A^\ell\right/{\mathcal I}d_A^{\ell+1}\right)(-2a-b)\to 0.\leqno{(3.3.3)_\ell}$$}

\noindent Since $a>0$ and $b\geq 0$, from the long exact cohomology sequence of $(3.3.3)_\ell$ and from $(3.3.2)$ we obtain
$$H^1\left(A, \left(\left.{\mathcal I}d_A^\ell\right/{\mathcal I}d_A^{\ell+1}\right)(V(-a))\right)=0\leqno{(3.3.4)}$$
for $\ell\geq 1$.  From the short exact sequence

{\footnotesize$$
0\to\left(\left.{\mathcal I}d_A^\ell\right/{\mathcal I}d_A^{\ell+1}\right)(V(-a))\to
\left(\left.{\mathcal O}_A\right/{\mathcal I}d_A^{\ell+1}\right)(V(-a))\to
\left(\left.{\mathcal O}_A\right/{\mathcal I}d_A^\ell\right)(V(-a))\to 0
$$}

\noindent it follows that

{\footnotesize$$
\Gamma\left(A,\left(\left.{\mathcal O}_A\right/{\mathcal I}d_A^{\ell+1}\right)(V(-a))\right)\to
\Gamma\left(A,\left(\left.{\mathcal O}_A\right/{\mathcal I}d_A^\ell\right)(V(-a))\right)\to
H^1\left(A,\left(\left.{\mathcal I}d_A^\ell\right/{\mathcal I}d_A^{\ell+1}\right)(V(-a))\right)$$}

\noindent is exact and
$$\Gamma\left(A,\left(\left.{\mathcal O}_A\right/{\mathcal I}d_A^{\ell+1}\right)(V(-a))\right)\to
\Gamma\left(A,\left(\left.{\mathcal O}_A\right/{\mathcal I}d_A^\ell\right)(V(-a))\right)
$$
is surjective for $\ell\geq 1$ by $(3.3.4)$.
Hence
$$
\Gamma\left(A,\left(\left.{\mathcal O}_A\right/{\mathcal I}d_A^\ell\right)(V(-a))\right)\to
\Gamma\left(A,{\mathcal O}_A(V(-a))\right)
$$
is surjective for $\ell\geq 2$.  This means that the element $\tau\in\Gamma\left(A,{\mathcal O}_A(V(-a))\right)$ can be extended
to a holomorphic section $\tilde\tau$ of $V(-a)$ on some open neighborhood $U$ of $A$ in ${\mathbb P}_n$.  By directly using the standard pseudoconvexity arguments for Hartogs's extension of holomorphic functions, we can now extend $\tilde\tau$ to a holomorphic section $\tau^\sharp\in\Gamma\left({\mathbb P}_n,V(-a)\right)$ of $V(-a)$ over ${\mathbb P}_n$.  Or, we can get the extension to $\tau^\sharp$ by using the result of Mathias Peternell that for any complex submanifold of complex dimension $n-q$ in ${\mathbb P}_n$ the complex manifold ${\mathbb P}_n-A$ is $(2q-1)$-complete [Peternell1987, p.433, Theorem] and by using [Andreotti-Grauert1962, p.254, Th\'eor\`eme 15].

\bigbreak\noindent(3.4) {\it Two Independent Sections.}  We claim that the holomorphic section $s\wedge\tau^\sharp$ of $(\det V)(-a)$ over ${\mathbb P}_n$ is not identically zero on ${\mathbb P}_n$, where $s$ is from (1.4).  Suppose the contrary.  Since the holomorphic section $s$ of $V$ over ${\mathbb P}_n$ is nowhere zero on ${\mathbb P}_n-Z_{\rm red}$, the identical vanishing of $s\wedge\tau^\sharp$ means that $\tau^\sharp=hs$ on ${\mathbb P}_n-Z_{\rm red}$ for some holomorphic section $h$ of ${\mathcal O}_{{\mathbb P}_n}(-a)$ over ${\mathbb P}_n-Z_{\rm red}$.  Since $Z_{\rm red}$ is of complex codimension $2$ in ${\mathbb P}_n$, it follows that the holomorphic section $h$ of ${\mathcal O}_{{\mathbb P}_n}(-a)$ over ${\mathbb P}_n-Z_{\rm red}$ can be extended to a holomorphic section $h^\sharp$ of ${\mathcal O}_{{\mathbb P}_n}(-a)$ over ${\mathbb P}_n$ and $\tau^\sharp=h^\sharp s$ holds everywhere on ${\mathbb P}_n$.  This implies that $\tau^\sharp$ is identically zero on $Z_{\rm red}$ and in particular identically zero on $A$, which is a contradiction.  An alternative argument to get a contradiction is that the positivity of the integer $a$ contradicts the existence of the non identically zero holomorphic section $h^\sharp$ of ${\mathcal O}_{{\mathbb P}_n}(-a)$ over ${\mathbb P}_n$.  Hence $s\wedge\tau^\sharp$ is a non identically zero holomorphic section of $(\det V)(-a)$ over ${\mathbb P}_n$, which contradicts
$(\det V)(-a)={\mathcal O}_{{\mathbb P}_n}(-b-a)$ and $a+b>0$.

\bigbreak

\bigbreak\noindent{\bf\S4. Branched Cover of Blowup of Complex Projective Space.}

\bigbreak\noindent(4.1) {\it General Case of No Regular Branch of Reduction of Codimension 2 Subspace.}  In the above argument for the special case of some branch $A=Z_{\rm red}^0$ of the reduction $Z_{\rm red}$ of the codimension 2 subspace $Z$ being regular, the regularity of $A=Z_{\rm red}^0$ is mainly used in the extension of the line bundle $[W]$ on $A$ to ${\mathbb P}_n$ and in the vanishing of $H^1\left(A, \left(\left.{\mathcal I}d_A^\ell\right/{\mathcal I}d_A^{\ell+1}\right)(-q)\right)$ for $\ell\geq 1$ and $q\geq 0$ by the theorem of Schneider-Zintl.  For the step of using the standard pseudoconvexity arguments for Hartogs's extension of holomorphic functions to extend $\tilde\tau$ to a holomorphic section $\tau^\sharp\in\Gamma\left({\mathbb P}_n,V(-a)\right)$ of $V(-a)$ over ${\mathbb P}_n$, the regularity of $A=Z_{\rm red}^0$ is not essential.  For the general case we have to find other ways to handle the two steps which depend on the regularity of $A=Z_{\rm red}^0$.  We are going to deal with the line bundle extension first.

\bigbreak\noindent(4.2) {\it Simply Connectedness of Normalization and Desingularization of Simply Connected Reduced Complex Space.}  Let $Y$ be a simply connected reduced complex space and let $\hat Y$ be the normalization of $Y$ and $\tilde Y$ be the desingularization of $\hat Y$ (by a finite number of successive monoidal transformations with nonsingular centers).  Then $\hat Y$ and $\tilde Y$ are simply connected.  

\medbreak The reason is that when we have a loop in either $\hat Y$ of $\tilde Y$, we can shrunk its image in $Y$ down to a point and we can do the shrinking by a map $\varphi$ to $Y$ from the closure of a domain $D$ with smooth boundary in ${\mathbb C}$ and assume after slightly perturbing the map $\varphi$ and changing the complex structure of $D$ that the map $\varphi$ maps $D$ holomorphically to $Y$ and the $\varphi$-image of $D$ is not entirely contained in the singular set of $Y$.  Any holomorphic map from $D\to Y$ can be uniquely lifted up to $\hat Y$ and $\tilde Y$, because the pullback of a weakly holomorphic function germ on $Y$ (which is a locally bounded function germ on $Y$ holomorphic on the regular part of $Y$) by $\varphi$ is holomorphic on $D$ and the lifting to $\hat Y$ can be obtained by using jets of $\varphi$ at points of $D$ of sufficiently high order.

\bigbreak\noindent(4.3) {\it Extension of Line Bundle from Singular Subvariety to Complex Projective Space.}  Mathias Peternell's Barth-Lefschetz theorem [Peternell1983, p.387, Satz 11] for singular subvarieties can be used to yield the following extension result.
Let $Y$ be an irreducible subvariety of complex dimension $n_Y$ in ${\mathbb P}_n$.
Let $\pi_{\hat Y}:\hat Y\to Y$ be the normalization of $Y$.
If $2\leq 2n_Y-n$, then for any holomorphic line bundle $L$ on $Y$ there exists $c\in{\mathbb Z}$ such that
the $\pi_{\hat Y}$-pullback of $L$ to $\hat Y$ is isomorphic to the the $\pi_{\hat Y}$-pullback of ${\mathcal O}_{{\mathbb P}_n}(c)$ to $\hat Y$.
In the application to our situation $Y$ is a branch $Z_{\rm red}^0$ of the reduction $Z_{\rm red}$ of the codimension 2 subspace $Z$ with $n_Y=n-2$ and $2\leq 2(n-2)-n=n-4$ or $n\geq 6$.

\medbreak The proof is as follows.  Since from $2\leq 2n_Y-n$ by the Barth-Lefschetz-Peternell theorem [Peternell1983, p.387, Satz 11] $H^2({\mathbb P}_n,{\mathbb Z})\to H^2(Y,{\mathbb Z})$ induced by the inclusion map $Y\hookrightarrow{\mathbb P}_n$ is surjective, there exists some $c\in{\mathbb Z}$ such that the restriction of ${\mathcal O}_{{\mathbb P}_n}(c)$ to $Y$ has the same Chern class as $L$.  By replacing $L$ by $L\otimes{\mathcal O}_{{\mathbb P}_n}(-c)$, we can assume without loss of generality that the Chern class of $L$ vanishes.  Since $Y$ is simply connected, its normalization $\hat Y$ and the desingularization $\tilde Y$ (obtained from $\hat Y$ by a finite number of successive monoidal transformations with nonsingular centers) are both simply connected.  Let $\hat L$ be the pullback of $L$ to $\hat Y$ and $\tilde L$ be the pullback of $\hat L$ to $\tilde Y$.  Since $\tilde Y$ is a compact K\"ahler manifold and the Chern class of $\tilde L$ vanishes, it follows that $\tilde L$ is a flat line bundle on $\tilde Y$.  Since $\tilde Y$ is simply connected, all flat bundles on $\tilde Y$ are trivial and in particular, $\tilde L$ is trivial and there is a global nowhere zero holomorphic section $F$ of $\tilde L$ over $\tilde Y$.  Let $\pi_{\tilde Y}:\tilde Y\to\hat Y$ be the desingularization projection map.  For any open subset $U$ of $\hat Y$ with $\hat L|_U$ trivial, the restriction of $F$ to $\left(\pi_{\tilde Y}\right)^{-1}(U)$ is a holomorphic function on $\left(\pi_{\tilde Y}\right)^{-1}(U)$ which must be the $\pi_{\tilde Y}$-pullback of a holomorphic function on $U$.  This means that $F$ is the $\pi_{\tilde Y}$-pullback of a nowhere zero holomorphic section of $\hat L$ over $\hat Y$.  Thus the line bundle $\hat L$ is trivial on $\hat Y$.

\bigbreak\noindent(4.4) {\it Vanishing Theorem for Singular Space.}  We now come to the question how to handle the difficulty coming from the second use of the regularity of $A=Z_{\rm red}^0$ in the vanishing of $H^1\left(A, \left(\left.{\mathcal I}d_A^\ell\right/{\mathcal I}d_A^{\ell+1}\right)(-q)\right)$ for $\ell\geq 1$ and $q\geq 0$ by the theorem of Schneider-Zintl in our argument in (3.3) for the special case of regular $A=Z_{\rm red}^0$.
For singular normal complex spaces there is the following vanishing theorem of Mumford, which is given on p.153 of the book of Shiffman-Sommese [Shiffman-Sommese1985] and which is proved by resolving the singularities.

\medbreak Let $L$ be a holomorphic line bundle over an irreducible normal projective algebraic variety $Y$.  Suppose there is a positive integer $m$ such that (i) there is no common zero in $Y$ for all elements of $H^0(Y, L^m)$ and (ii) the complex dimension of the image of the holomorphic map $Y\to{\mathbb P}_N$ whose components are a ${\mathbb C}$-basis of $H^0(Y, L^m)$ is $k>1$. Then $H^p\left(Y, L^*\right)=0$ for $p<\min(2,k)$.

\bigbreak As a corollary we have the following vanishing theorem for use in our situation.   Let $Y$ be an irreducible (reduced) complex space of complex dimension $\geq 2$ and $\pi:\hat Y\to Y$ be the normalization of $Y$.  Let $L$ be a holomorphic line bundle on $Y$ which is ample.  Then $H^1\left(\hat Y,\pi^*L\right)=0$.

\bigbreak\noindent(4.5) {\it Branched Cover over Blowup of ${\mathbb P}_n$ along Complex Line.}
Since both the line bundle extension and the vanishing theorem for the case of a singular $A=Z_{\rm red}^0$ requires the use of its normalization, we are going to construct now a branched cover $\hat{\mathcal E}$ over the blow-up ${\mathcal E}$ of ${\mathbb P}_n$ along a complex line in ${\mathbb P}_n$ so that the normalization $\hat A$ of $A$ is the inverse image of $A$ under the branched-cover map $\hat{\mathcal E}\to{\mathcal E}$ composed with the blow-up map ${\mathcal E}\to{\mathbb P}_n$.  First we construct ${\mathcal E}$ as follows.

\bigbreak\noindent(4.6) {\it Light-Source Projection.}  For homogeneous coordinates $\left[\zeta_0,\cdots,\zeta_n\right]$ of ${\mathbb P}_n$ denote by ${\mathbb P}_{n-1}^\infty$ the hyperplane $\left\{\zeta_0=0\right\}$ at infinity.  Denote by  $L^\infty_{1,2}$ the complex line $\left\{\zeta_0=\zeta_3=\zeta_4=\cdots=\zeta_n=0\right\}$ inside ${\mathbb P}_{n-1}^\infty$.  It consists of all infinity points of the $(x_1,x_2)$-coordinate plane in the affine coordinates $x_1=\frac{\zeta_1}{\zeta_0},\cdots,x_n=\frac{\zeta_n}{\zeta_0}$ of the affine part ${\mathbb C}^n={\mathbb P}_n-{\mathbb P}_{n-1}^\infty$ of ${\mathbb P}_n$.  We choose the 
$(x_1,x_2)$-coordinate plane simply because it is the most natural one to use.

\medbreak Now we choose a system of homogeneous coordinates $\left[\zeta_0,\cdots,\zeta_n\right]$ of ${\mathbb P}_n$ such that ${\mathbb P}_{n-1}^\infty\cap Z_{\rm red}$ is of complex codimension $2$ in ${\mathbb P}_{n-1}^\infty$ and disjoint from $L^\infty_{1,2}$.
For $P\in{\mathbb P}_n$ not in $L^\infty_{1,2}$ let $\Pi_P$ be the complex plane in ${\mathbb P}_n$ containing $P$ and $L^\infty_{1,2}$.
Let $T_{1,2}$ be the complex projective space $\left\{\zeta_1=\zeta_2=0\right\}$ of complex dimension $n-2$ in inside ${\mathbb P}_n$, which is disjoint from $L^\infty_{1,2}$.

\medbreak Define the {\it light-source projection} with light source $L^\infty_{1,2}$ and target $T_{1,2}$ as the map $\pi_{1,2}:{\mathbb P}_n-L^\infty_{1,2}\to T_{1,2}$ which sends $P$ to the intersection point $\pi_{1,2}(P)$ of $T_{1,2}$ and $\Pi_P$.
The fibers of $\pi_{1,2}$ are $\Pi_P\cap\left({\mathbb P}_n-{\mathbb P}_{n-1}^\infty\right)$, which are biholomorphic to ${\mathbb C}^2$ so that $\pi_{1,2}:{\mathbb P}_n-L^\infty_{1,2}\to T_{1,2}$ is a holomorphic plane bundle over $T_{1,2}$.
We add $L^\infty_{1,2}$ to each fiber $\Pi_P\cap\left({\mathbb P}_n-{\mathbb P}_{n-1}^\infty\right)$ of $\pi_{1,2}:{\mathbb P}_n-L^\infty_{1,2}\to T_{1,2}$ to get a ${\mathbb P}_2$-bundle $\pi_{\mathcal E}: {\mathcal E}\to T_{1,2}$ with fibers $\Pi_P$.  The ${\mathbb P}_2$-bundle $\pi_{\mathcal E}: {\mathcal E}\to T_{1,2}$
is the same as the space obtained from ${\mathbb P}_n$ by blowing up the complex line $L^\infty_{1,2}$.

\bigbreak\noindent(4.7) {\it Construction of Branched Cover over Blowup of ${\mathbb P}_n$.}  The light-source map $\pi_{1,2}:{\mathbb P}_n-L^\infty_{1,2}\to T_{1,2}$ makes the branch $A=Z^0_{\rm red}$ a branched cover over $T_{1,2}$ because $Z_{\rm red}$ is disjoint from $L^\infty_{1,2}$.  Let $\hat\pi:\hat A \to A$ be the normalization of $A=Z^0_{\rm red}$.
For every point $\hat P$ of $\hat A$, take the complex plane $\Pi_{\hat\pi(\hat P)}$ which contains $\hat\pi(\hat P)$ and $L^\infty_{1,2}$ and use it as the fiber and put them together as $\hat P$ varies in $\hat A$ to form the ${\mathbb P}_2$-bundle $\pi_{\hat{\mathcal E}}: \hat{\mathcal E}\to \hat A$.
We use the identity map of $\Pi_P$ to map each fiber of $\hat{\mathcal E}$ to the complex plane $\Pi_P$ of ${\mathbb P}_n$ to get a branched-cover map $\Pi_{\hat{\mathcal E}}:\hat{\mathcal E}\to {\mathcal E}$.   The map $\Pi_{\hat{\mathcal E}}$ restricted to $\hat A$ is the normalization map $\hat A\to A$.

\bigbreak\noindent(4.8) {\it Alternative Description of Branched Cover over Blowup of ${\mathbb P}_n$.}   The branched cover $\hat{\mathcal E}$ of ${\mathcal E}$ is biholomorphic to the pullback of the ${\mathbb P}_2$-bundle $\pi_{\mathcal E}: {\mathcal E}\to T_{1,2}$ over $T_{1,2}$ by the branched-over map $\pi_{1,2}\circ\hat\pi: \hat A\to T_{1,2}$.  This can be seen by considering how the complex planes $\Pi_{\hat\pi(\hat P)}$ are put together to form $\hat{\mathcal E}$ as $\hat P$ varies in $\hat A$ and how the complex planes $\Pi_P$ are put together to form ${\mathcal E}$ as $P$ varies in $T_{1,2}$.

\medbreak When the branched cover $\hat{\mathcal E}$ of ${\mathcal E}$ is regarded as the pullback of the ${\mathbb P}_2$-bundle $\pi_{\mathcal E}: {\mathcal E}\to T_{1,2}$ over $T_{1,2}$ by the branched-over map $\pi_{1,2}\circ\hat\pi: \hat A\to T_{1,2}$, there is a natural projection $\Pi^\prime_{\hat{\mathcal E}}: \hat{\mathcal E}\to{\mathcal E}$ whose restriction to the subvariety $\hat A$ of $\hat{\mathcal E}$ agrees with the map $\pi_{1,2}\circ\hat\pi: \hat A\to T_{1,2}$ (when $T_{1,2}$ is regarded as a subvariety of ${\mathcal E}$).  We would like to point out that this natural projection $\Pi^\prime_{\hat{\mathcal E}}: \hat{\mathcal E}\to{\mathcal E}$ is different from the branched cover map $\Pi_{\hat{\mathcal E}}:\hat{\mathcal E}\to {\mathcal E}$ in (4.7).  The difference is the following.  For a point $P$ of ${\mathbb P}_n$ not in $L^\infty_{1,2}$, when the identification of $\Pi_P$ with the fiber $\pi_{\hat{\mathcal E}}^{-1}(\hat P)$ of $\hat{\mathcal E}$ with $\hat\pi(\hat P)\in A\cap\Pi_P$ and the identification of $\Pi_P$ with the fiber $\pi_{\mathcal E}^{-1}(\Pi_P\cap T_{1,2})$ of ${\mathcal E}$ are so chosen that the restriction of $\Pi_{\hat{\mathcal E}}:\hat{\mathcal E}\to {\mathcal E}$ to the fiber $\pi_{\hat{\mathcal E}}^{-1}(\hat P)$ of $\hat{\mathcal E}$ is the identity map of $\Pi_P$, then the restriction of the map $\Pi^\prime_{\hat{\mathcal E}}:\hat{\mathcal E}\to {\mathcal E}$ to the fiber $\pi_{\hat{\mathcal E}}^{-1}(\hat P)$ of $\hat{\mathcal E}$ induces the map on $\Pi_P$ which is the identity map on $\Pi_P\cap L^\infty_{1,2}$ but is a translation of $\Pi_P-L^\infty_{1,2}$ in the affine coordinates of $\Pi_P-L^\infty_{1,2}$ mapping $\hat\pi(\hat P)$ to the only point of $\Pi_P\cap T_{1,2}$.

\bigbreak\noindent(4.9) {\it Pullbacks to Branched Cover over Blowup of ${\mathbb P}_n$ along Complex Line.}  We use $\pi_{\hat{\mathcal E}}\circ\pi_{\mathcal E}:\hat{\mathcal E}\to{\mathbb P}_n$ to pull back $V$ to $\hat V$, to pull back $s\in\Gamma({\mathbb P}_n,V)$ to $\hat s\in\Gamma(\hat{\mathcal E},\hat V)$, to pull back $\sigma^0\in\Gamma(A, V)$ to $\hat\sigma^0\in\Gamma(\hat A,\hat V)$, and to pull back the divisor $W$ of $\sigma^0$ in $A=Z_{\rm red}^0$ to the divisor $\hat W$ of $\hat\sigma^0$ on the normalization $\hat A$ of $A$.

\medbreak By the above result on extending line bundles from singular subvarieties to ${\mathbb P}_n$ in (4.3), there exists some integer $a$ such that on $\hat A$ the $\hat\pi$-pullback of the line bundle $[W]$ on $A$ agrees with the $\hat\pi$-pullback of the line bundle ${\mathcal O}_{{\mathbb P}_n}(a)$ on ${\mathbb P}_n$.  Since the intersection of $A$ and the hypersurface $H$ of ${\mathbb P}_n$ is contained in $W$ according to (3.1), it follows that $a>0$.  Let $t_{\hat W}$ be the canonical section of the line bundle $[\hat W]$ on $\hat A$.  Let $\hat\tau=\frac{\hat\sigma^0}{t_{\hat W}}$.  Then $\hat\tau$ is a holomorphic section of $\hat V(-a)$ over $\hat A$ which is nowhere zero on $A$.  Here $\hat V(-a)$ is obtained by twisting $\hat V$ by the $\hat\pi$-pullback of the line bundle ${\mathcal O}_{{\mathbb P}_n}(-a)$ on ${\mathbb P}_n$.

\medbreak The important advantage of introducing the branched cover $\hat{\mathcal E}$ of ${\mathcal E}$ is that though both the complex space $\hat{\mathcal E}$ and its subvariety $\hat A$ are singular in general, the normal bundle $N_{\hat A,\hat{\mathcal E}}$ of $\hat A$ in $\hat{\mathcal E}$ is isomorphic to the pullback of the normal bundle of $N_{T_{1,2},\,{\mathbb P}_n}$ of $T_{1,2}$ in ${\mathbb P}_n$ under $\pi_{1,2}\circ\hat\pi: \hat A\to T_{1,2}$ so that $N_{\hat A, \hat{\mathcal E}}$ is the direct sum of two positive line bundles.  By Mumford's vanishing theorem for negative line bundles on normal complex spaces, we have the vanishing $$H^1\left(\hat A, {\rm Sym}^\ell((N_{\hat A, \hat{\mathcal E}})^*)(-q)\right)=0$$ for $\ell\geq 1$ and $q\geq 0$ (when $n\geq 5$). Here again ${\rm Sym}^\ell((N_{\hat A, \hat{\mathcal E}})^*)(-q)$ is obtained by twisting it by the $\hat\pi$-pullback of the the line bundle ${\mathcal O}_{{\mathbb P}_n}(-q)$ on ${\mathbb P}_n$.  Let ${\mathcal I}d_{\hat A}$ be the full ideal sheaf of $\hat A$ on $\hat{\mathcal E}$.  Since $\left({\mathcal I}d_{\hat A}\right)^\ell\left/\left({\mathcal I}d_{\hat A}\right)^{\ell+1}\right.$ is isomorphic to the sheaf of germs of holomorphic sections of ${\rm Sym}^\ell((N_{\hat A, \hat{\mathcal E}})^*)$, it follows that
$$H^1\left(\hat A, \left(\left({\mathcal I}d_{\hat A}\right)^\ell\left/\left({\mathcal I}d_{\hat A}\right)^{\ell+1}\right.\right)(-q)\right)=0\leqno{(4.9.1)}$$ for $\ell\geq 1$ and $q\geq 0$.
At this point we can just use straightforward modifications to repeat as follows the above arguments for the special case of $A=Z_{\rm red}^0$ being regular given in \S3.

\bigbreak\noindent(4.10) {\it Adaption of Final Argument to General Case.}  We have a short exact sequence
$$
0\to{\mathcal O}_{\hat A}\stackrel{\hat\Phi}{\longrightarrow}{\mathcal O}_{\hat A}(\hat V(-a))\to{\mathcal O}_{\hat A}(-2a-b)\to 0\leqno{(4.10.1)}
$$ on $\hat A$, where the inclusion map $\hat\Phi$ is defined by multiplication by $\hat\tau$.  Tensoring the short exact sequence $(4.10.1)$ with ${\rm Sym}^\ell\left(\left(N_{\hat A,\hat{\mathcal E}}\right)^*\right)$, we get the short exact sequence

{\footnotesize$$0\to\left.{\mathcal I}d_{\hat A}^\ell\right/{\mathcal I}d_{\hat A}^{\ell+1}\to\left(\left.{\mathcal I}d_{\hat A}^\ell\right/{\mathcal I}d_{\hat A}^{\ell+1}\right)(\hat V(-a))\to\left(\left.{\mathcal I}d_{\hat A}^\ell\right/{\mathcal I}d_{\hat A}^{\ell+1}\right)(-2a-b)\to 0\leqno{(4.10.2)_\ell}$$}

\noindent on $\hat A$.  Since $a>0$ and $b\geq 0$, from the long exact cohomology sequence of $(4.10.2)_\ell$ and from $(4.9.1)$ we obtain
$$H^1\left(\hat A, \left(\left.{\mathcal I}d_{\hat A}^\ell\right/{\mathcal I}d_{\hat A}^{\ell+1}\right)(\hat V(-a))\right)=0\leqno{(4.10.3)}$$
for $\ell\geq 1$.  From the short exact sequence

{\footnotesize$$
0\to\left(\left.{\mathcal I}d_{\hat A}^\ell\right/{\mathcal I}d_{\hat A}^{\ell+1}\right)(\hat V(-a))\to
\left(\left.{\mathcal O}_{\hat A}\right/{\mathcal I}d_{\hat A}^{\ell+1}\right)(\hat V(-a))\to
\left(\left.{\mathcal O}_{\hat A}\right/{\mathcal I}d_{\hat A}^\ell\right)(\hat V(-a))\to 0
$$}

\noindent on $\hat A$ it follows that

{\footnotesize$$
\Gamma\left(\hat A,\left(\left.{\mathcal O}_{\hat A}\right/{\mathcal I}d_{\hat A}^{\ell+1}\right)(\hat V(-a))\right)\to
\Gamma\left(\hat A,\left(\left.{\mathcal O}_{\hat A}\right/{\mathcal I}d_{\hat A}^\ell\right)(\hat V(-a))\right)\to
H^1\left(\hat A,\left(\left.{\mathcal I}d_{\hat A}^\ell\right/{\mathcal I}d_{\hat A}^{\ell+1}\right)(\hat V(-a))\right)$$}

\noindent is exact and
$$\Gamma\left(\hat A,\left(\left.{\mathcal O}_{\hat A}\right/{\mathcal I}d_{\hat A}^{\ell+1}\right)(\hat V(-a))\right)\to
\Gamma\left(\hat A,\left(\left.{\mathcal O}_{\hat A}\right/{\mathcal I}d_{\hat A}^\ell\right)(\hat V(-a))\right)
$$
is surjective for $\ell\geq 1$ by $(4.10.3)$.
Hence
$$
\Gamma\left(\hat A,\left(\left.{\mathcal O}_{\hat A}\right/{\mathcal I}d_{\hat A}^\ell\right)(\hat V(-a))\right)\to
\Gamma\left(\hat A,{\mathcal O}_{\hat A}(\hat V(-a))\right)
$$
is surjective for $\ell\geq 2$.  This means that the element $\hat\tau\in\Gamma\left(\hat A,{\mathcal O}_{\hat A}(\hat V(-a))\right)$ can be extended
to a holomorphic section $\widetilde{\hat{\tau}}$ of $\hat V(-a)$ on some connected open neighborhood $\hat U$ of $\hat A$ in $\hat{\mathcal E}$.  Since the subvariety $A$ in ${\mathbb P}_n$ is disjoint from $L^\infty_{1,2}$, the complement $\hat A$ of $\left(\pi_{\hat{\mathcal E}}\circ\pi_{\mathcal E}\right)^{-1}\left(L^\infty_{1,2}\right)$ in $\hat{\mathcal E}$ is a ${\mathbb C}^2$-bundle $\hat{\mathcal E}^0$ over $\hat A$ whose associated ${\mathbb P}_2$-bundle (as its compactification) is the ${\mathbb P}_2$-bundle $\hat{\mathcal E}$ over $\hat A$.  The ${\mathbb C}^2$-bundle $\hat{\mathcal E}^0$ over $\hat A$ is the direct sum of two positive line bundles over $\hat A$.  The subvariety $\hat A$ in $\hat{\mathcal E}$ is the zero-section of the ${\mathbb C}^2$-bundle $\hat{\mathcal E}^0$ over $\hat A$ with $\hat U$ being an open neighborhood of its zero-section $\hat A$.  By the standard argument of extending holomorphic functions locally over any strict $(n-2)$-pseudoconcave boundary in an $n$-dimensional normal complex space, we can extend the holomorphic section $\widetilde{\hat{\tau}}$ of $\hat V(-a)$ from $\hat U$ to a holomorphic section $\hat\tau^\natural$ of $\hat V(-a)$ on $\hat{\mathcal E}^0$.

\medbreak We claim that the holomorphic section $\hat s\wedge\hat\tau^\natural$ of $(\det\hat V)(-a)$ over $\hat{\mathcal E}^0$ is not identically zero on $\hat{\mathcal E}^0$.  Suppose the contrary.  Since the holomorphic section $\hat s$ of $\hat V$ over $\hat{\mathcal E}^0$ is nowhere zero on $\left(\pi_{\hat{\mathcal E}}\circ\pi_{\mathcal E}\right)^{-1}\left({\mathbb P}_n-\left(Z_{\rm red}\cup L^\infty_{1,2}\right)\right)$, the identical vanishing of $\hat s\wedge\hat\tau^\natural$ means that $\hat\tau^\natural=\hat h \hat s$ on $\left(\pi_{\hat{\mathcal E}}\circ\pi_{\mathcal E}\right)^{-1}\left({\mathbb P}^n-\left(Z_{\rm red}\cup L^\infty_{1,2}\right)\right)$ for some holomorphic section $\hat h$ of ${\mathcal O}_{\hat{\mathcal E}}(-a)$ over $\left(\pi_{\hat{\mathcal E}}\circ\pi_{\mathcal E}\right)^{-1}\left({\mathbb P}^n-\left(Z_{\rm red}\cup L^\infty_{1,2}\right)\right)$.  Since $\left(\pi_{\hat{\mathcal E}}\circ\pi_{\mathcal E}\right)^{-1}\left(Z_{\rm red}\right)$ is of complex codimension $2$ in $\hat{\mathcal E}^0$, it follows that the holomorphic section $\hat h$ of ${\mathcal O}_{\hat{\mathcal E}}(-a)$ over $\left(\pi_{\hat{\mathcal E}}\circ\pi_{\mathcal E}\right)^{-1}\left({\mathbb P}^n-\left(Z_{\rm red}\cup L^\infty_{1,2}\right)\right)$ can be extended to a holomorphic section $\hat h^\natural$ of $\hat h$ of ${\mathcal O}_{\hat{\mathcal E}}(-a)$ over $\hat{\mathcal E}^0$ and $\hat\tau^\natural=\hat h^\natural \hat s$ holds everywhere on $\hat{\mathcal E}^0$.  This implies that $\hat\tau^\natural$ is identically zero on $\left(\pi_{\hat{\mathcal E}}\circ\pi_{\mathcal E}\right)^{-1}\left(Z_{\rm red}\right)$ and in particular identically zero on $\hat A$, which is a contradiction.  Hence $\hat s\wedge\hat\tau^\natural$ is a non identically zero holomorphic section of $(\det \hat V)(-a)$ over $\hat{\mathcal E}^0$.

\medbreak Let $\lambda$ be the number of sheets in the the branched-cover map $\Pi_{\hat{\mathcal E}}:\hat{\mathcal E}\to {\mathcal E}$.  The restriction of $\Pi_{\hat{\mathcal E}}$ to $\hat{\mathcal E}^0$ is a branched-cover map of $\hat{\mathcal E}^0$ over ${\mathbb P}_n-L^\infty_{1,2}$ of $\lambda$ sheets.  The line bundle $(\det \hat V)(-a)$ over $\hat{\mathcal E}$ is the pullback of the line bundle $(\det V)(-a)$ on ${\mathbb P}_n$ by the map $\pi_{\hat{\mathcal E}}\circ\pi_{\mathcal E}:\hat{\mathcal E}\to{\mathbb P}_n$.  For a generic point $P$ of ${\mathbb P}_n-L^\infty_{1,2}$, define the element $\xi(P)$ of the fiber of $\left((\det V)(-a)\right)^{\otimes\lambda}$ over $P$ by
$$
\xi(P)=\prod\left\{\left(\hat s\wedge\hat\tau^\natural\right)(\hat P)\,\Big|\, \hat P\in\hat{\mathcal E}^0\ \ {\rm with}\ \left(\pi_{\hat{\mathcal E}}\circ\pi_{\mathcal E}\right)(\hat P)=P\,\right\}
$$
to obtain a non identically zero section $\xi$ of $\left((\det V)(-a)\right)^{\otimes\lambda}$ over ${\mathbb P}_n-L^\infty_{1,2}$.  In other words, $\xi\in\Gamma\left({\mathbb P}_n-L^\infty_{1,2},\left((\det V)(-a)\right)^{\otimes\lambda}\right)$ is the multiplicative direct image of $\hat s\wedge\hat\tau^\natural\in\Gamma\left(\hat{\mathcal E}^0, (\det \hat V)(-a)\right)$ under the map $\pi_{\hat{\mathcal E}}\circ\pi_{\mathcal E}$.
Since $L^\infty_{1,2}$ is of codimension $n-1\geq 2$ in ${\mathbb P}_n$, we can extend $$\xi\in H^0\left({\mathbb P}_n-L^\infty_{1,2}, \left((\det V)(-a)\right)^{\otimes\lambda})\right)$$ to a non identically zero holomorphic section of $\left((\det V)(-a)\right)^{\otimes\lambda}$ over ${\mathbb P}_n$, which contradicts $\left((\det V)(-a)\right)^{\otimes\lambda}={\mathcal O}_{{\mathbb P}_n}\left(\lambda(-b-a)\right)$ on ${\mathbb P}_n$ with  $a+b>0$.  This finishes the proof of the Main Theorem.

\bigbreak\noindent{\bf References}

\bigbreak\noindent[Andreotti-Grauert1962] Aldo Andreotti and Hans Grauert,
Th\'eor\`eme de finitude pour la cohomologie des espaces complexes.
{\it Bull. Soc. Math. France} \textbf{90} (1962), 193 –- 259.

\medbreak\noindent
[Barth1970] Wolf Barth,
Transplanting cohomology classes in complex-projective space.
{\it Amer. J. Math.} \textbf{92} (1970), 951 –- 967.

\medbreak\noindent
[Barth1975] Wolf Barth,
Larsen's theorem on the homotopy groups of projective manifolds of small embedding codimension. {\it Algebraic geometry} (Proc. Sympos. Pure Math., Vol. 29, Humboldt State Univ., Arcata, Calif., 1974), pp. 307 –- 313. Amer. Math. Soc., Providence, R.I., 1975.

\medbreak\noindent
[Barth1979] Wolf Barth, Review of the paper of Hans Grauert and Michael Schneider,
``Komplexe Unterräume und holomorphe Vektorraumb\"undel vom Rang zwei'' [{\it Math. Ann.} \textbf{230} (1977), 75 –- 90].
{\it Mathematical Reviews} \textbf{58} (1979), \#1279 (MR0481134).

\medbreak\noindent
[Barth-Larsen1972] Wolf Barth and Mogens Esrom Larsen,
On the homotopy groups of complex projective algebraic manifolds.
{\it Math. Scand.} \textbf{30} (1972), 88 –- 94.

\medbreak\noindent
[Faltings1981] Gerd Faltings,
Verschwindungss\"atze und Untermannigfaltigkeiten kleiner Kodimension des komplex-projektiven Raumes.
{\it J. Reine Angew. Math.} \textbf{326} (1981), 136 –- 151.

\medbreak\noindent
[Grauert-Schneider1977] Hans Grauert and Michael Schneider,
Komplexe Unterr\"aume und holomorphe Vektorraumb\"undel vom Rang zwei.
{\it Math. Ann.} \textbf{230} (1977),  75 -– 90.

\medbreak\noindent
[Larsen1973] Mogens Esrom Larsen,
On the topology of complex projective manifolds.
{\it Invent. Math.} \textbf{19} (1973), 251 –- 260.

\medbreak\noindent
[LePotier1975] Joseph Le Potier,
Annulation de la cohomolgie \`a valeurs dans un fibr\'e vectoriel holomorphe positif de rang quelconque.
{\it Math. Ann.} \textbf{218} (1975), 35 –- 53.

\medbreak\noindent[Nakano1955] Shigeo Nakano,
On complex analytic vector bundles.
{\it J. Math. Soc. Japan} \textbf{7} (1955), 1 –- 12.

\medbreak\noindent[Peternell1983] Mathias Peternell,
Ein Lefschetz-Satz f\"ur Schnitte in projektiv-algebraischen Mannigfaltigkeiten.
{\it Math. Ann.} \textbf{264} (1983),  361 -– 388.

\medbreak\noindent[Peternell1987] Mathias Peternell,
$q$-completeness of subsets in complex projective space.
{\it Math. Z.} \textbf{195} (1987), 443 –- 450.

\medbreak\noindent
[Okonek-Schneider-Spindler2011] Christian Okonek, Michael Schneider, and Heinz Spindler,
{\it Vector bundles on complex projective spaces} (with Appendix by S. I. Gelfand). Birkh\"auser/Springer 2011.

\medbreak\noindent
[Schneider1975] Michael Schneider,
Lefschetzs\"atze und Hyperkonvexit\"at.
{\it Invent. Math.} \textbf{31} (1975),  183 –- 192.

\medbreak\noindent
[Schneider1987] Michael Schneider,
Vector bundles and submanifolds of projective space: nine open problems. Algebraic geometry, Bowdoin, 1985 (Brunswick, Maine, 1985), 101 –- 107,
{\it Proc. Sympos. Pure Math.} \textbf{46}, Part 2, Amer. Math. Soc., Providence, RI, 1987.

\medbreak\noindent
[Schneider-Zintl1993] Michael Schneider and J\"org Zintl,
The theorem of Barth-Lefschetz as a consequence of Le Potier's vanishing theorem.
{\it Manuscripta Math.} \textbf{80} (1993), 259 –- 263.

\medbreak\noindent[Shiffman-Sommese1985]
Bernard Shiffman and Andrew John Sommese,
{\it Vanishing theorems on complex manifolds}.
Progress in Mathematics, \textbf{56}. Birkh\"auser Boston, Inc., Boston, MA, 1985.

\medbreak\noindent[Sommese1979]
Andrew John Sommese,
Complex subspaces of homogeneous complex manifolds. I. Transplanting theorems.
{\it Duke Math. J.} \textbf{46} (1979), 527 –- 548.

\medbreak\noindent[Sommese1982]
Andrew John Sommese,
Complex subspaces of homogeneous complex manifolds. II. Homotopy results.
{\it Nagoya Math. J.} \textbf{86} (1982), 101 –- 129.

\medbreak\noindent[vandeVen1980] Antonius van de Ven, Twenty years of classifying algebraic vector bundles, In: {\it Algebraic Geometry
Angers 1979}, Sijthoff and Noordhoff International Publishers (1980), 3 -- 20.

\bigbreak\noindent{\it Author's address:} Department of Mathematics, Harvard University, Cambridge, MA 02138, U.S.A.

\bigbreak\noindent{\it Author's e-mail address:} siu@math.harvard.edu

\end{document}